\documentclass{article}
\usepackage[arrow, matrix, curve]{xy}
\usepackage{amsthm}
\usepackage{hyperref}
\input{amssym.tex}
\newtheorem{thm}{Theorem}

\newtheorem{obs} {Observation}
\newtheorem{pro} {Proposition}

\newtheorem{ass} {Assumption}

	\newenvironment{pf}{{\it Proof:}\quad}{\hfill$QED$}

\begin{document}

	\title{Mutation and recombination for hyperbolic 3-manifolds}
	\author{Thilo Kuessner}
	\date{}
	\maketitle
	\begin{abstract}  We give a short topological proof for Ruberman's Theorem about mutation and volume, using the Maskit combination theorem and the homology of the linear group.
%	\footnotetext{2000 Mathematics Subject Classification }
\noindent
		\end{abstract}

Let $\rho:\Gamma\rightarrow SL\left(2,{\bf C}\right)$ be a faithful representation with discrete, torsion-free image, 
and let $M=\Gamma\backslash{\bf H}^3$ be the associated orientable, hyperbolic 3-manifold. We assume that $\Gamma$ has finite covolume, that is $vol\left(M\right)<\infty$. 

Let $\Sigma\subset M$ be a properly embedded, incompressible, boundary-incompressible, 2-sided surface. Let $\tau:\Sigma\rightarrow \Sigma$ be a diffeomorphism. We will denote 
by $M^\tau$ the result of cutting $M$ along $\Sigma$ and regluing via $\tau$. 

Ruberman's Theorem (\cite[Theorem 1.3]{rub}) states that for some specific pairs $\left(\Sigma,\tau\right)$, for example for the hyperelliptic involution of the genus 2 surface, one always has that $M^\tau$ is hyperbolic and $vol\left(M^\tau\right)=vol\left(M\right)$.

The pairs $\left(\Sigma,\tau\right)$ considered in \cite{rub} have in common that $\tau$ has finite order and that for each $\rho$ as above 
the representation $\rho\mid_{\pi_1\Sigma}\circ \tau_*:\pi_1\Sigma\rightarrow SL\left(2,{\bf C}\right)$ is conjugate in $SL\left(2,{\bf C}\right)$ to $\rho\mid_{\pi_1\Sigma}$, see \cite[Theorem 2.2]{rub}.

As indicated in \cite[Theorem 4.4]{rub} this is actually the only specific property of the pairs $\left(\Sigma,\tau\right)$ which is needed. That is, the conclusion of $M^\tau$ being hyperbolic and having the same volume as $M$ does actually hold for all pairs $\left(\Sigma,\tau\right)$ such that $\rho\mid_{\pi_1\Sigma}\circ \tau_*$ is conjugate in $SL\left(2,{\bf C}\right)$ to $\rho\mid_{\pi_1\Sigma}$ (and such that $\Sigma$ is not a virtual fiber).

The proof in \cite{rub} uses hard results about minimal surfaces. As Ruberman points out in \cite[Theorem 4.4]{rub}, an alternative proof is possible using Bers' Theorem about the deformation space of quasifuchsian groups, yet another hard analytical theorem. 

The aim of this paper is to give a purely topological proof of Ruberman's result which does not use any facts from analysis, but rather a simple version of the Maskit combination theorems and group-homological arguments. (The combination theorems have a purely topological proof in \cite{mas}, although nowadays there exist also
analytical proofs using minimal surfaces or harmonic maps.)\\
\\
{\bf Theorem 1}: {\em Let $M$ be a compact, orientable, connected $3$-manifold with (possibly empty) boundary. Let $\rho:\Gamma\rightarrow SL\left(2,{\bf C}\right)$ be a lift of a faithful
representation $\Gamma\rightarrow Isom^+\left({\bf H}^3\right)=PSL
\left(2,{\bf C}\right)$ with discrete, torsion-free image. Assume that
$\rho\left(\Gamma\right)\backslash{\bf H}^3$ has finite volume and is diffeomorphic to $int\left(M\right)=M-\partial M$.

Let $\Sigma\subset M$ be a properly embedded, connected, incompressible, boundary-incompressible, 2-sided surface which is not a virtual fiber.

Let $\tau:\left(\Sigma,\partial\Sigma\right)\rightarrow \left(\Sigma,\partial\Sigma\right)$ be an orientation-preserving diffeomorphism of pairs such that $\tau^m=id$ for some $m\in{\bf N}$ and
such that there
exists some $A\in SL\left(2,{\bf C}\right)$ with $\rho\left(\tau_*h\right)=A\rho\left(h\right)A^{-1}$
for all $h\in\pi_1\Sigma$.

Then $int\left(M^\tau\right)$ is hyperbolic and $vol\left(M^\tau\right)=vol\left(M\right)$.}\\

The proof of $vol\left(M^\tau\right)=vol\left(M\right)$ will follow an argument analogous to the proof of \cite[Theorem 2.13]{neu} 
but will use the constructions from \cite[Section 4]{ku} to handle the case of manifolds with several cusps (where relative group homology does not directly apply).

\section{Preliminaries}

Let $M$ be an orientable, complete, hyperbolic 3-manifold of finite volume and $\Sigma$ a properly embedded, 2-sided surface. Let $\Gamma=\pi_1M$. (We assume without further mention $M$ and $\Sigma$ to be connected and we omit basepoints from the notation.) The monodromy has image in $Isom^+\left({\bf H}^3\right)=PSL\left(2,{\bf C}\right)$.
By Culler's Theorem \cite[Corollary 2.2]{cul} the monodromy comes from 
a representation $\rho:\Gamma\rightarrow SL\left(2,{\bf C}\right)$ and we will henceforth assume such a representation $\rho$ to be fixed. 

Let $\tau:\Sigma\rightarrow \Sigma$ be a diffeomorphism. Throughout the paper we will work with the following assumption on $\left(\Sigma,\tau\right)$:

\begin{ass}\label{ass} There exists some $A\in SL\left(2,{\bf C}\right)$ with $$\rho\left(\tau_*h\right)=A\rho\left(h\right)A^{-1}$$
for all $h\in\pi_1\Sigma$.\end{ass}

An obvious example where this condition is satisfied, is the case that $\Sigma$ is totally geodesic and $\tau$ is an isometry of the induced metric. In this case one can upon conjugation assume that $\rho\left(\pi_1\Sigma\right)\subset SL\left(2,{\bf R}\right)$ and then $A$ is an elliptic element in $SL\left(2,{\bf R}\right)\subset SL\left(2,{\bf C}\right)$. 

Less obvious examples are provided by \cite[Theorem 2.2]{rub} which asserts that \hyperref[ass]{Assumption \ref*{ass}} holds whenever $\left(\Sigma,\tau\right)$ is what \cite{rub} calls a symmetric surface, for example if $\Sigma$ has genus $2$ and $\tau$ is an hyperelliptic involution, or for certain symmetries of the 3- or 4-punctured sphere or the 1- or 2-punctured torus. 
For these symmetric surfaces $\left(\Sigma,\tau\right)$, \hyperref[ass]{Assumption \ref*{ass}} holds 
regardless how $\Sigma$ is embedded (as an incompressible, boundary-incompressible surface) into a hyperbolic $3$-manifold $M$.\\

Let $M^\tau$ be the result of cutting $M$ along $\Sigma$ and regluing via $\tau$. Since $\Sigma$ is a 2-sided, properly embedded surface, it has a 
neighborhood $N\simeq\Sigma\times\left[0,1\right]$ in $M$, and a neighborhood $N^\tau\simeq\Sigma\times\left[0,1\right]$ in $M^\tau$. The complements $M-int\left(N\right)$ and $M^\tau-int\left(N^\tau\right)$ are diffeomorphic and we let $X$ be
the union of $M$ and $M^\tau$ along this identification of $M-int\left(N\right)$ and $M^\tau-int\left(N^\tau\right)$. The union of $N$ and $N^\tau$ yields a copy of the mapping torus $T^\tau$ in $X$.

Let $\pi_1M=<S\mid R>$ be a presentation of $\pi_1M$.
The Seifert-van Kampen Theorem implies
$$\pi_1X=<S,t\mid R, tht^{-1}=\tau_*\left(h\right)\ \forall\ h\in\pi_1\Sigma>.$$

If \hyperref[ass]{Assumption \ref*{ass}} holds, then we have a well-defined
representation $$\rho_X:\pi_1X\rightarrow SL\left(2,{\Bbb C}\right)$$
by $\rho_X\left(s\right)=\rho\left(s\right)$ for all $s\in S$ and $\rho_X\left(t\right)=A$.

Composition with the homomorphisms $\pi_1M\rightarrow \pi_1X$ and $\pi_1
M^\tau\rightarrow \pi_1X$ induced by the inclusions yields the given representation $\rho$
of $\pi_1M$ and a representation $\rho^\tau$ of $\pi_1M^\tau$. We will show in the Section 2 that under
certain hypotheses (namely that $\Sigma$ is incompressible, boundary-incompressible and not a virtual fiber) the representation $\rho^\tau$ 
will be discrete and faithful and therefore $M^\tau$ is a hyperbolic manifold (although $X$ is not).\\
%\footnotemark\footnotetext[1]{However $X$ clearly is not a hyperbolic manifold. In fact $A$ is of finite order, therefore $\rho_X\left(\pi_1X\right)\subset SL\left(2,{\bf C}\right)$ is not torsion-free.}.

For use in Section 2 we mention that the Seifert-van Kampen Theorem also allows a description of $\pi_1M^\tau$ (as a subgroup of $\pi_1X$ with the above presentation), where we have to distinguish the cases whether the surface $\Sigma$ is a separating surface in $M$ or not. 

If an incompressible, 2-sided
surface $\Sigma$ separates $M$ into two submanifolds $M_1$ and $M_2$, then $\pi_1M=\pi_1M_1*_{\pi_1\Sigma}\pi_1M_2$ is the amalgamated product of $\pi_1M_1$ and $\pi_1M_2$ over $\pi_1\Sigma$ for two 
monomorphisms $\phi_1:\pi_1\Sigma\rightarrow \pi_1M_1,\phi_2:\pi_1\Sigma\rightarrow \pi_1M_2$. For the mutation $M^\tau$ we obtain that
$\pi_1M^\tau$ 
is the amalgamated product of $\pi_1M_1$ and $t\pi_1M_2t^{-1}$ amalgamated over $\pi_1\Sigma$ via the monomorphisms $\psi_1:\pi_1\Sigma\rightarrow \pi_1M_1, \psi_2:\pi_1\Sigma\rightarrow \pi_1M_2$ defined by
$$\psi_1\left(h\right)=\phi_1\left(\tau_*h\right), \psi_2\left(h\right)=
t\phi_2\left(h\right)t^{-1}$$
for all $h\in \pi_1\Sigma$.

If an incompressible, 2-sided surface $\Sigma$ does not separate $M$, then there is a manifold $N$ 
with $\partial N=\Sigma_1\cup\Sigma_2$ such that $M$ is the quotient of $N$ under 
some diffeomorphism $f:\Sigma_1\rightarrow\Sigma_2$. We assume the basepoint $x_0$ to belong to $\Sigma_1$, let 
$u$ be some path in $N$ from $x_0$ to $f\left(x_0\right)$, let $K_1=\pi_1\Sigma_1$
and let $K_2\subset \pi_1N$ the subgroup of all based homotopy classes of loops of the form $u*\sigma*\overline{u}$, where 
$\sigma:\left[0,1\right]\rightarrow\Sigma_2$ with $\sigma\left(0\right)=\sigma\left(1\right)=f\left(x_0\right)$. Then $\pi_1M=\pi_1N*_\alpha$ 
is the HNN-extension of $\pi_1N$ for the two monomorphisms $\phi_1:K_1\rightarrow \pi_1N,\phi_2:K_2\rightarrow \pi_1N$ induced by the inclusions 
and the isomorphism $\alpha:K_1\rightarrow K_2$ defined by $\alpha\left(\left[\sigma\right]\right)=
\left[u*f\left(\sigma\right)*\overline{u}\right]$, where $\left[\sigma\right]$ is
the based homotopy class of $\sigma:\left[0,1\right]\rightarrow\Sigma_1$ with $\sigma\left(0\right)=\sigma\left(1\right)=x_0$. Let $v$ be the extending element. For the
mutation $M^\tau$ we obtain that
$\pi_1M^\tau=\pi_1N*_{\alpha\tau_*}$ is the HNN-extension of $\pi_1N$, with extending element $vt$, for the same monomorphisms $\phi_1,\phi_2$ and the isomorphism $\alpha\tau_*:K_1\rightarrow K_2$.\\

We mention that $X$ clearly
is not a manifold, let alone a hyperbolic one. In fact, the representation $\rho_X$ is not faithful and its image is not torsion-free, as the following elementary observation shows.

\begin{obs}\label{elliptic} Let $\rho:\Gamma\rightarrow SL\left(2,{\bf C}\right)$ be a discrete, faithful representation of a torsion-free group, let $M=\rho\left(\Gamma\right)\backslash{\bf H}^3$ be the associated hyperbolic 3-manifold, and let $\Sigma\subset M$ be a properly embedded, incompressible,
boundary-incompressible surface. Assume that $\Sigma$ is not boundary-parallel.

Let $\tau:\Sigma\rightarrow\Sigma$ be a diffeomorphism such that $\tau^m=id$ for some $m\in{\bf N}$ and such that
there exists some $A\in SL\left(2,{\bf C}\right)$ with $\rho\left(\tau_*h\right)=A\rho\left(h\right)A^{-1}$
for all $h\in\pi_1\Sigma$.

Then $A^m=\pm{\bf 1}$.\end{obs}
\begin{pf} We have $A^m\rho\left(h\right)A^{-m}=A^{m-1}\rho\left(\tau_*h\right)A^{-\left(m-1\right)}=\ldots=\rho\left(\tau_*^mh\right)=\rho\left(h\right)$ for all $h\in\pi_1\Sigma$, thus $A^m$ conjugates $\rho\left(\pi_1\Sigma\right)$ to itself.

If $\Sigma$ is not boundary-parallel, then the image of $\pi_1\Sigma$ in $Isom^+\left({\bf H}^3\right)$ contains at least two non-commuting hyperbolic elements. The corresponding matrices in $\rho\left(\pi_1\Sigma\right)\subset SL\left(2,{\bf C}\right)$ are then diagonalisable
with respect to two distinct bases. 

A matrix that conjugates a diagonalizable 2-by-2-matrix (with 2 distinct eigenvalues) to itself must be diagonalizable with respect to the same basis. Thus
the conjugating matrix $A^m\in SL\left(2,{\bf C}\right)$ must be diagonal with respect to two distinct bases.
This implies, by elementary linear algebra, that $A^m$ is a multiple of the identity. Together with $det\left(A^m\right)=1$ this enforces $A^m=\pm{\bf 1}$.\end{pf}\\
\\
%Remark: If \hyperref[ass]{Assumption \ref*{ass}} holds for a matrix $A\in SL\left(2,{\bf C}\right)$, then it also holds for each multiple of $A$ by some $\xi\in{S}^1

%The idea to recover Rubermans result via a fundamental class construction seems to appear for the first time in \cite[Theorem 2.13]{neu}, where Neumann also considers the effect of mutation on images of integer fundamental classes. (In the context of hyperbolic 3-manifolds, where one can use known computations of Chern-Simons-invariants.)

\section{Discreteness of recombinations}

Let $M$ be an
oriented (and connected) compact 3-manifold with boundary such that its interior $int\left(M\right)=M-\partial M$ is hyperbolic. 
The hyperbolic structure is given by the conjugacy class of a representation $\pi_1M\rightarrow PSL\left(2,{\Bbb C}\right)$ which by Culler's Theorem \cite[Corollary 2.2]{cul} is the composition of a representation $\rho:\pi_1M\rightarrow SL\left(2,{\Bbb C}\right)$ with the canonical projection $PSL\left(2,{\Bbb C}\right)\rightarrow SL\left(2,{\Bbb C}\right)$. For a diffeomorphism $\tau:\Sigma\rightarrow\Sigma$ we let $\rho^\tau:\pi_1M^\tau
\rightarrow SL\left(2,{\Bbb C}\right)$ be the representation defined in Section 1.

\begin{pro}\label{discrete} If $\Sigma\subset M$ is a
properly embedded, connected, incompressible, boundary-incompressible, 2-sided surface, which is
not a virtual fiber, and if \hyperref[ass]{Assumption \ref*{ass}} holds for a diffeomeorphism $\tau:\Sigma\rightarrow\Sigma$, then $\rho^\tau\left(\pi_1M^\tau\right)$ is a discrete subgroup of $SL\left(2,{\Bbb C}\right)$. In particular, $int\left(M^\tau\right)$ is hyperbolic.\end{pro}

\begin{pf} Since $\Sigma$ is incompressible, boundary-incompressible
and
not a virtual fiber, $H:=\rho\left(\pi_1\Sigma\right)$ is geometrically finite by the Thurston-Bonahon Theorem\footnotemark\footnotetext[1]{The Thurston-Bonahon Theorem is a combination of results in \cite{mar},\cite{thu},\cite{bon}. An explicit statement can be found in \cite[Theorem 1.1]{coo} or \cite[Corollary 8.3]{can}.}. This implies that 
the limit set $W=\Lambda\left(H\right)\subset \partial_\infty{\bf H}^3$ is a quasi-circle, in particular a Jordan curve. The Sch\"onflies theorem implies that $W$ decomposes $ \partial_\infty{\bf H}^3\cong {\bf S}^2$ into two topological disks.

We distinguish the cases that $\Sigma$ is separating in $M$ (and thus in $M^\tau$) or not. We wish to apply the Maskit Combination Theorems \cite[Chapter VII]{mas} and have to check that their assumptions hold for $\rho^\tau\left(\pi_1M^\tau\right)$.

{\bf Case 1: $\Sigma$ is separating.} Then we obtain from Section 1 that $G=\rho\left(\pi_1M\right)$ is an 
amalgamated 
product $G=G_1*_HG_2$
with $G_1=\rho\left(\pi_1M_1\right)$ and $G_2=\rho\left(\pi_1M_2\right)$. 
The complement of the limit set $\partial_\infty{\bf H}^3-W$ consists of two open topological disks $B_1,B_2$ such that $B_i$ is precisely $H$-invariant in $G_i$ for $i=1,2$, that is
$$h\left(B_i\right)\subset B_i \ \ \forall\ h\in H$$
$$g\left(B_i\right)\subset B_{3-i} \ \ \forall\ g\in G_i - H.$$

To apply the first Maskit combination theorem we have to check that the analogous condition is satisfied for the decomposition of
$G^\tau:=\rho^\tau\left(\pi_1M^\tau\right)$ as amalgamated product of $G_1$ and $AG_2A^{-1}$ over $H$.

We already know that $h\left(B_i\right)\subset B_i$ for $h\in H$ and $g\left(B_1\right)\subset B_{2}$ for $g\in G_1-H$. We have to check that $AgA^{-1}\left(B_2\right)\subset B_1$ for $g\in G_2-H$.

\hyperref[ass]{Assumption \ref*{ass}} implies that conjugation by $A\in SL\left(2,{\bf C}\right)$ maps $H$ to $H$, thus $A$ and $A^{-1}$ map the limit set $W=\Lambda\left(H\right)$ to itself. Since $A$ is orientation-preserving it must then also map $B_1$ to $B_1$ and $B_2$ to $B_2$. Hence $g\left(B_2\right)\subset B_1$ implies $AgA^{-1}\left(B_2\right)
\subset B_1$ for $g\in G_2-H$.

The first combination theorem implies then that $G^\tau$ is discrete.

{\bf Case 2: $\Sigma$ is non-separating.}
Then we obtain from Section 1 that $G=\rho\left(\pi_1M\right)$ is an HNN-extension $G=G_0*_\alpha$ of $G_0:=\rho\left(\pi_1N\right)$
for some isomorphism $\alpha:H\rightarrow H_2$, where $H:=\rho\left(K_1\right)=\rho\left(\pi_1\Sigma_1\right), H_2:=\rho\left(K_2\right)$. 

Let $f=\rho\left(v\right)\in G$ be the extending element of the HNN-extension $G=G_0*_\alpha$ and let $W_2=f\left(W\right)\subset\partial_\infty{\bf H}^3$.
Then $W$ and $W_2$ are disjoint simple closed curves and the complement 
$\partial_\infty{\bf H}^3-\left(W\cup W_2\right)$ consists of two open topological disks $B_1,B_2$ and an open annulus $R$ such that $B_i$ is precisely $H_i$-invariant in $G_0$ for $i=1,2$ and, moreover, 
$$f\left(R\cup B_2\right)= B_2.$$

To apply the second Maskit combination theorem we have to check that the analogous conditions are satisfied for the description of 
$G^\tau:=\rho^\tau\left(\pi_1M^\tau\right)$ as an HNN-extension $G_\tau=G_0*_{\alpha\tau_*}$ for the isomorphism $\alpha\tau_*:H\rightarrow H_2$, with extending element $\rho\left(vt\right)=fA$.

This is clear for the first condition: we already know that $B_i$ is precisely $H_i$-invariant in $G_0$.

\hyperref[ass]{Assumption \ref*{ass}} implies again that $A$ and $A^{-1}$ map the limit set $W=\Lambda\left(H\right)$ to itself. Since $A$ is orientation-preserving it must then also map $B_1$ bijectively to $B_1$ and $R\cup B_2$ bijectively to $R\cup B_2$. In particular $f\left(R\cup B_2\right)= B_2$ implies
$$fA\left(R\cup B_2\right)= B_2.$$
%Moreover $\tau_*$ and $\tau_*^{-1}$ map $W_1$ into $W_1$
%and $W_2$ into $W_2$ and hence (since it is orientation-preserving) also $B_1$ into $B_1$, $B_2$ into $B_2$, $R$ into $R$. Thus $$\alpha\tau_*\left(R\cup B_2\right)\subset B_2, \left(\alpha\tau_*\right)^{-1}\left(R\cup B_1\right)=\tau_*^{-1}\alpha^{-1}\left(R\cup B_1\right)\subset B_1.$$

The second combination theorem implies then that $G^\tau$ is discrete.

\end{pf}\\

One may wonder whether it is possible to apply the topological combination theorems (\cite[Theorem VII.A.12.,C.13.]{mas}) directly to the actions on ${\bf H}^3$ rather than to the actions on $\partial_\infty{\bf H}^3$. Letting $\tilde{\phi}:
{\bf H}^3\rightarrow {\bf H}^3$ be the "lift" of $\tau$ defined as the lift of the isometry $\phi:\rho\left(\pi_1\Sigma\right)\backslash{\bf H}^3\rightarrow 
\rho\left(\pi_1\Sigma\right)\backslash{\bf H}^3$ from \cite[Lemma 2.5.]{rub} this would 
require to have a $\tilde{\phi}$-invariant copy of $\widetilde{\Sigma}$ in ${\bf H}^3$. The main part of the argument in \cite{rub} actually
consists in providing such a $\tilde{\phi}$-invariant surface, using heavy machinery from the theory of minimal surfaces. The advantage of using the combination theorem for the action on $\partial_\infty{\bf H}^3$ instead of the action on ${\bf H}^3$
is that one can avoid this machinery.

\section{Volume of mutations}
\subsection{Recollections}
We recollect some constructions from \cite[Section 4]{ku} which will be needed in the proof of Theorem 1, in particular to handle the 
case of disconnected boundary, where relative group homology does not directly apply. 

Let $\left(X,Y\right)$ be a pair of topological spaces. Assume that $X$ is path-connected and $Y$ has path-components 
$Y_1,\ldots,Y_s$. Let $\Gamma=\pi_1\left(X,x\right)$. Using pathes from $y_i$ to $x$ one can fix isomorphisms 
$l_i:\pi_1\left(Y_i,y_i\right)
\rightarrow \Gamma_i$ to subgroups $\Gamma_i\subset\Gamma$ for $i=1,\ldots,s$.

As in \cite[Section 4.2.1]{ku}, we will denote by $DCone\left(\cup_{i=1}^sY_i\rightarrow X\right)$ the union along $Y$
of $X$ with the {\em disjoint} cones over $Y_1,\ldots,Y_s$. Let $\Psi_X:X\rightarrow \mid B\Gamma\mid$ be the classifying 
map for $\pi_1X$. (We consider $B\Gamma$ as a simplicial set and denote $\mid B\Gamma\mid$ its geometric realization.)

As in \cite[Definition 7]{ku} we denote $B\Gamma^{comp}=DCone\left(\cup_{i=1}^s B\Gamma_i\rightarrow B\Gamma\right)$, that is the
quasi-simplicial set which is the union along $\cup_{i=1}^sB\Gamma_i$ of the simplicial set $B\Gamma$ and the cones $Cone\left(B\Gamma_i\right)$ over $B\Gamma_i$ with cone point $c_i$ for $i=1,\ldots,s$.
$\Psi_X$ extends to a continuous map $\Psi_X:DCone\left(\cup_{i=1}^sY_i\rightarrow X\right)\rightarrow \mid B\Gamma^{comp}\mid$. 
If $X$ and $Y_1,\ldots,Y_s$ are aspherical, then composition of $\left(\Psi_X\right)_*$ with the isomorphism $\iota:H_*\left( \mid B\Gamma^{comp}\mid\right)\rightarrow H_*^{simp}\left(B\Gamma^{comp}\right)$ yields the inverse of the Eilenberg-MacLane isomorphism:
$$\iota\circ \left(\Psi_X\right)_*=EM^{-1}:H_*\left(DCone\left(\cup_{i=1}^sY_i\rightarrow X\right)\right)\rightarrow H_*^{simp}\left(B\Gamma^{comp}\right).$$
In particular, if $M$ is a compact, connected, oriented, aspherical  manifold with aspherical boundary $\partial M=\partial_1M\cup\ldots\cup\partial_sM$, then there
is for $*\ge 2$
a canonical isomorphism\footnotemark\footnotetext[2]{An explicit realisation of this isomorphism is 
as follows. Let $z\in Z_*\left(M,\partial M\right)$ be a relative cycle, then 
$\partial z\in Z_{*-1}\left(\partial M\right)$ and thus $z+Cone\left(\partial z\right)\in Z_*\left(DCone\left(\cup_{i=1}^s\partial_iM\rightarrow M\right)\right)$ is a cycle. This induces an isomorphism in homology.} $H_*\left(M,\partial M\right)\cong H_*\left(DCone\left(\cup_{i=1}^s\partial_iM\rightarrow M\right)\right)$  
and we obtain the isomorphism (\cite[Lemma 8]{ku})
$$EM^{-1}:H_*\left(M,\partial M\right)\rightarrow H_*^{simp}\left(B\Gamma^{comp}\right).$$

Now assume that $int\left(M\right)=M-\partial M$ is hyperbolic and 
$\rho:\Gamma\rightarrow SL\left(2,{\bf C}\right)$ 
is a representation defining the hyperbolic structure. As 
in \cite[Definition 6]{ku} let 
$$BSL\left(2,{\bf C}\right)^{comp}=
DCone\left(\dot{\cup}_{c\in\partial_\infty{\bf H}^3}BSL\left(2,{\bf C}\right)\rightarrow BSL\left(2,{\bf C}\right)\right),$$ where each 
of the cone points of the disjoint copies of $Cone\left(BSL\left(2,{\bf C}\right)\right)$ is identified with a point $c\in\partial_\infty{\bf H}^3$.

Assume that $vol\left(M\right)<\infty$. Then for each $i=1,\ldots,s$ there is a unique $c_i\in\partial_\infty{\bf H}^3$ with $\Gamma_i\subset Fix\left(c_i\right)$. Therefore $B\rho:B\Gamma\rightarrow BSL\left(2,{\bf C}\right)$ can be extended to
$$B\rho:B\Gamma^{comp}\rightarrow BSL\left(2,{\bf C}\right)^{comp}$$
by mapping the cone point of $Cone\left(B\Gamma_i\right)$ to the 
cone point $c_i\in\partial_\infty{\bf H}^3$ of $BSL\left(2,{\bf C}\right)^{comp}$ with $\Gamma_i\subset Fix\left(c_i\right)$. 

%In the case of hyperbolic $3$-manifolds, the Eilenberg-MacLane isomorphism $EM$ can be explicitly defined 
%(for a fixed $\tilde{x}_0\in{\bf H}^3$) by the chain map $\Psi:C_*^{simp}\left(B\Gamma^{comp}\right)\rightarrow
%C_*\left(DCone\left(\cup_{i=1}^s\partial_iM\rightarrow M\right)\right)$ by
%$$\Psi\left(g_1,\ldots,g_k\right)=\pi\left(str\left(\tilde{x}_0,g_1\tilde{x}_0,\ldots,g_1\ldots g_k\tilde{x}_0\right)\right)\mbox{\ \ if\ \ }\left(g_1,\ldots,g_k\right)\in C_*^{simp}\left(B\Gamma\right)$$
%$$\Psi\left(p_1,\ldots,p_{k-1},c_i\right)=\pi\left(str\left(\tilde{x}_0,p_1\tilde{x}_0,\ldots,p_1\ldots p_{k-1}\tilde{x}_0,c_i\right)\right)\mbox{\ \ if\ \ }\left(p_1,\ldots,p_{k-1},c_i\right)\in C_*^{simp}\left(Cone\left(B\Gamma_i\right)\right),$$
%see \cite[Lemma 8a]{ku}.
%(The condition $\Gamma_i\subset Fix\left(c_i\right)$ is needed here to guarantee that $\Psi$ is a chain map.)

%$$\left(B\rho\right)_*EM^{-1}\left[M,\partial M\right]\in H_3^{simp}\left(BSL\left(2,{\bf C}\right)^{comp};{\bf Q}\right).$$

In \cite[Section 4.2.3]{ku} we defined a cocycle $\overline{c\nu_3}
\in C_{simp}^3\left(BSL\left(2,{\bf C}\right)^{comp};{\bf R}\right)$ whose cohomology class $\left[\overline{c\nu_3} \right]$, by
the proof of \cite[Theorem 4]{ku}, satisfies $$<\left[\overline{c\nu_3}\right], \left(B\rho\right)_*EM^{-1}\left[M,\partial M\right]>=vol\left(M\right).$$

In particular, $vol\left(M\right)$ is determined by the element $$\left(B\rho\right)_*EM^{-1}\left[M,\partial M\right]\in H_3^{simp}\left(BSL\left(2,{\bf C}\right)^{comp};{\bf Q}\right).$$

\subsection{Proof of Theorem 1}
\begin{thm}\label{cusped}
Let $M$ be a compact, orientable, connected $3$-manifold with boundary $\partial M$. Let $\rho:\Gamma\rightarrow SL\left(2,{\bf C}\right)$ be a lift of a faithful 
representation $\Gamma\rightarrow Isom^+\left({\bf H}^3\right)=PSL
\left(2,{\bf C}\right)$ with discrete, torsion-free image. Assume that 
%$\rho\left(\Gamma\right)$ acts properly discontinuously on ${\bf H}^3$ and that 
$\rho\left(\Gamma\right)\backslash{\bf H}^3$ has finite volume and is diffeomorphic to $int\left(M\right)=M-\partial M$.

Let $\Sigma\subset M$ be a properly embedded,
connected, incompressible, boundary-incompressible, 2-sided surface which is not a virtual fiber.

Let $\tau:\left(\Sigma,\partial\Sigma\right)\rightarrow \left(\Sigma,\partial\Sigma\right)$ be an orientation-preserving diffeomorphism of pairs such that $\tau^m=id$ for some $m\in{\bf N}$ and
such that there
exists some $A\in SL\left(2,{\bf C}\right)$ with $\rho\left(\tau_*h\right)=A\rho\left(h\right)A^{-1}$
for all $h\in\pi_1\Sigma$.

Then $int\left(M^\tau\right)$ is hyperbolic and $vol\left(M^\tau\right)=vol\left(M\right)$.
\end{thm}

\begin{pf}
We have proved in \hyperref[discrete]{Proposition \ref*{discrete}} that $M^\tau$ is hyperbolic, that is the representation $\rho^\tau:\pi_1M^\tau\rightarrow SL\left(2,{\bf C}\right)$ is faithful and has discrete image.

Let $X$ be constructed as in Section 1 and let $\partial X=X\cap\left(\partial M\cup\partial M^\tau\right)$. Let $\Gamma^X:=\pi_1X$ and let $\rho_X:\Gamma^X
\rightarrow SL\left(2,{\bf C}\right)$ be the representation, defined in Section 1, with $\rho_X\mid_{\Gamma}=\rho$ and $\rho_X\left(t\right)=A$.

The construction of $X$ implies that
\begin{equation}\label{a}i_{M*}\left[M,\partial M\right]-i_{M^\tau *}\left[M^\tau,\partial M^\tau\right]= i_{T^\tau *}\left[T^\tau,\partial T^\tau\right]\in H_3\left(X,\partial X;{\bf Q}\right),\end{equation}
where $i_M,i_{M^\tau},i_{T^\tau}$ are the inclusions of $M,M^\tau$ and the mapping torus $T^\tau$ into $X$, and $\left[M,\partial M\right],\left[M^\tau,\partial M^\tau\right],\left[T^\tau,\partial T^\tau\right]$ are the fundamental classes in homology with ${\bf Q}$-coefficients.\\
%The hyperbolic manifolds $M,M^\tau$ and the mapping torus $T^\tau$ are Eilenberg-MacLane spaces, thus we get (unique up to homotopy) continuous maps
%$B\rho:M\rightarrow BSL\left(2,{\bf C}\right), B\rho^\tau:M^\tau\rightarrow BSL\left(2,{\bf C}\right), B\rho_{T}:T^\tau\rightarrow BSL\left(2,{\bf C}\right)$.
%Recall from Section 1 that \hyperref[ass]{Assumption \ref*{ass}}
%yields a representation $\rho_X: \pi_1X\rightarrow SL\left(2,{\bf C}\right)$
%whose restrictions to $\pi_1M,\pi_1M^\tau,\pi_1T^\tau$ give $\rho,\rho^\tau,\rho_T$.
%Let $\Phi_X:X\rightarrow BSL\left(2,{\bf C}\right)$ be the composition of the classifying map $X\rightarrow B\pi_1X$ with $B\rho_X:B\pi_1X\rightarrow BSL\left(2,{\bf C}\right)$.
%Then we have $\Phi_X\circ i_M\sim B\rho$ and $\Phi_X\circ i_{M^\tau}\sim B\rho^\tau$ because the induced homomorphisms of the fundamental groups agree.
%Thus we have $$B\rho_*\left[M,\partial M\right]-B\rho^\tau_*\left[M^\tau,\partial M^\tau\right]=B\rho_{X*}\left(i_{M*}\left[M\right]\right)-B\rho_{X*}\left(i_{M^\tau *}\left[M^\tau\right]\right)$$
%$$
%=B\rho_{X*}\left(i_{T^\tau *}\left[T^\tau,\partial T^\tau\right]\right)=B\rho_{T*}\left[T^\tau,\partial T^\tau\right].$$

We observe that the path-components $\partial_1 M,\ldots,\partial_sM$ of $\partial M$ are in 1-1-correspondence with the path-components $\partial_1M^\tau,\ldots,\partial_sM^\tau$ of $\partial M^\tau$ 
and with the path-components $\partial_1X,\ldots,\partial_sX$ of $\partial X$. To each path-component $\partial_j\Sigma$ of $\partial \Sigma$ there is some path-component $\partial_{i_j}M$ of $\partial M$ with $\partial_j\Sigma\subset \partial_{i_j}M$. 
The path-components $\partial_j\Sigma$ are in 1-1-correspondence with the path-components $\partial_jT^\tau$ of $T^\tau$, and we have then $\partial_j T^\tau\subset\partial_{i_j}X$. 

We choose the base point $x$ of $\Gamma=\pi_1\left(M,x\right)$ to belong to $M\cap M^\tau\cap T^\tau$, for $i=1,
\ldots,s$ the base point $x_i$ of $\pi_1\left(\partial_iM,x_i\right)$ to 
belong to $\partial_iM\cap\partial_iM^\tau\cap\partial T^\tau$, and for all $j$
the base point $y_j$ of $\pi_1\left(\partial_j T^\tau,y_j\right)$ to belong to $\partial_jT^\tau \cap \partial_{i_j}M\cap\partial_{i_j}M^\tau$.

Let $i\in\left\{1,\ldots,s\right\}$. Recall that a path $p_i$ from $x_i$ to $x$ yields an isomorphism $l_i:\pi_1\left(X,x_i\right)\rightarrow \pi_1\left(X,x\right)$. Let $\Gamma_i^X=l_i\left(\pi_1\left(\partial_iX,x_i\right)\right)\subset \Gamma^X:=\pi_1\left(X,x\right)$, that is $\Gamma_i^X$ is an isomorphic image of $
\pi_1\left(\partial_iX,x_i\right)$ in $\Gamma^X$. 
We can choose the path $p_i$ in $M\cap M^\tau$, thus if $i_{\partial_iM*}:\pi_1\left(\partial_iM,x_i\right)\rightarrow
\pi_1\left(\partial_iX,x_i\right)$ is induced by the inclusion $i_{\partial_iM}:\partial_iM\rightarrow \partial_iX$, then 
$j_{\partial_iM}l_i\mid_{\pi_1\left(\partial_iM,x_i\right)}=l_ii_{\partial_iM *}$, where 
$j_{\partial_iM}:\Gamma_i\rightarrow\Gamma_i^X$ is the inclusion.
%We can choose the path $p_i$ in $M\cap M^\tau$, thus $l_i$ maps $\pi_1\left(\partial_iM,x_i\right)$ to a subgroup $\Gamma_i$ of $\Gamma$ resp.\ $\pi_1\left(\partial_iM^\tau,x_i\right)$ to a subgroup $\Gamma_i^\tau$
%of $\Gamma^\tau:=\pi_1\left(M^\tau,x\right)$. That is we have $j_{\partial_iM}l_i=l_ii_{\partial_iM *}$, where $i_{\partial_iM}:\partial_iM\rightarrow\partial_iX$ and $j_{\partial_iM}:\Gamma_i\rightarrow\Gamma_i^X$ is the inclusion. 
An analogous fact holds by replacing $M$ with $M^\tau$.
%and then use the same path  . Thus, if we use $l_i$ to identify $\pi_1\left(\partial_iM,x_i\right)$ with a subgroup of $\Gamma=\pi_1\left(M,x\right)$ resp.\ 
%to identify $\pi_1\left(\partial_iM^\tau,x_i\right)$ with a subgroup of $\Gamma^\tau:=\pi_1\left(M^\tau,x\right)$, then we have $\Gamma_i\subset\Gamma_i^X$ and $\Gamma_i^\tau\subset\Gamma_i^X$.

For each $y_j\in \partial_jT^\tau\cap\partial_{i_j}M\cap\partial_{i_j}M^\tau$ we can choose a path 
from $y_j$ to $x_{i_j}$
in $\partial_{i_j}M\cap\partial_{i_j}M^\tau$. This yields an isomorphism $l_j$ from $\pi_1\left(\partial_jT^\tau,y_j\right)$ to a subgroup of 
$\pi_1\left(\partial_{i_j}M,x_{i_j}\right)$. Composition of $l_j$ with $l_{i_j}$ yields an isomorphism
from $\pi_1\left(\partial_jT^\tau,y_j\right)$ to a subgroup $\Gamma_j^T$ of $\Gamma$. The same pathes can be used to construct an isomorphism 
from $\pi_1\left(\partial_jT^\tau,y_j\right)$ to a subgroup of $\Gamma^\tau$, and both images are then
by construction the same subgroup of $\Gamma^X$.\\
% , $y\in T^\tau\cap M\cap M^\tau$ and a path from $y$ to $y_j$ in $T^\tau\cap M\cap M^\tau$ to get an isomorphism from $\pi_1\partial_jT^\tau$ w%ith a subgroup $\Gamma_j^T$ of $\Gamma^T:=\pi_1\left(T^\tau,y\right)$.\\

%Thus, if we denote $j_M,j_{M^\tau},j_{T^\tau}$ the homomorphisms of fundamental groups induced by $i_M,i_{M^\tau},i_{T^\tau}$, then $j_M$ maps $\Gamma_i$ to $\Gamma_i^X$ and $j_{M^\tau}$ maps $\Gamma_i^\tau$ to $\Gamma_i^X$. Hence we obtain well-defined maps $Bj_M:B\Gamma^{comp}\rightarrow \left(B\Gamma^X\right)^{comp}$ and $Bj_{M^\tau}:\left(B\Gamma^\tau\right)^{comp}\rightarrow\left(B\Gamma^X\right)^{comp}$.
%Application of $EM^{-1}$ yields, by naturality of the Eilenberg-MacLane isomorphism, that
%\begin{equation}\label{b}\left(Bj_M\right)_*EM^{-1}\left[M,\partial M\right]-\left(Bj_{M^\tau}\right)_*EM^{-1}\left[M^\tau,\partial M^\tau\right]=
%\left(Bj_{T^\tau}\right)_*EM^{-1}\left[T^\tau,\partial T^\tau\right]\end{equation}
%in $H_3\left(\left(B\Gamma^X\right)^{comp};{\bf Q}\right)$.\\

By Section 3.1 $B\rho$ 
%and $B\rho^\tau$ 
extends to a simplicial map $B\rho:B\Gamma^{comp}\rightarrow BSL\left(2,{\bf C}\right)^{comp}$ 
%and $
%B\rho:B\Gamma^{comp}\rightarrow BSL\left(2,{\bf C}\right)^{comp}$ 
where the cone point over $B\Gamma_i$ 
%resp.\ $B\Gamma_i^\tau$ are 
is mapped to the (unique) $c_i$ with $\Gamma_i\subset Fix\left(c_i\right)$.
% and $Fix\left(\Gamma_i^\tau\right)\subset c_i$.
However $B\rho_X:B\Gamma^X\rightarrow BSL\left(2,{\bf C}\right)$ does not extend to $\left(B\Gamma^X\right)^{comp}$ in such 
a way that restriction
to $B\Gamma^{comp}$ would give back $B\rho$. Indeed in general there is no $c_i\in \partial_\infty{\bf H}^3$ with 
$\rho_X\left(\Gamma_i^X\right)\subset Fix\left(c_i\right)$. 
(This is because $\rho_X\left(t\right)=A$ is not parabolic but elliptic.) Therefore we have to\footnotemark\footnotetext[3]{If $M$ is a closed manifold, then one can use $B\rho_X$ as in
\cite[Theorem 2.13]{neu} and does not need to go to a finite cover.} 
go to a finite cover as follows.

Let $\pi_1X=<S,t\mid R, tht^{-1}=\tau_*\left(h\right)\ \forall\ h\in\pi_1\Sigma>$ be the presentation from Section 1.
% and $\rho_X:\pi_1X\rightarrow SL\left(2,{\bf C}\right)$ the representation with $\rho_X\left(t\right)=A$.  By \hyperref[elliptic]{Observation \ref*{elliptic}} we have $A^m={\bf 1}$. Upon replacing the matrix $A$ in \hyperref[ass]{Assumption \ref*{ass}} by $iA$ (if $m$ even) resp.\ $-A$ (if $m$ odd) 
%we can assume that $A^m={\bf 1}$ and that \hyperref[ass]{Assumption \ref*{ass}} still holds.
Then $$a\left(
t\right)=1, a\left(s\right)=0\ \forall\ s\in S$$ yields a well-defined, surjective homomorphism
$$a:\pi_1X\rightarrow {\bf Z}/2m{\bf Z}.$$ Let $\pi:\widehat{X}\rightarrow X$ 
be the $2m$-fold cyclic covering with $\Gamma^{\widehat{X}}:=\pi_1\left(\widehat{X},\hat{x}\right)\cong ker\left(a\right)$ for some $\hat{x}\in\pi^{-1}\left(x\right)$. 
For each $i\in\left\{1,\ldots,s\right\}$ we let $\partial_i\widehat{X}\subset\widehat{X}$ be the preimage of $\partial_iX$, we fix 
the preimage $\hat{p}_i$ of the path $p_i$ ending at $\hat{x}$, and use $\hat{p}_i$ to define the isomorphism from $\pi_1\partial_i\widehat{X}$ to a subgroup
$\Gamma_i^{\widehat{X}}\subset \Gamma^{\widehat{X}}$. Then $\pi\left(\hat{p}_i\right)=p_i$ implies $\pi_*\left(\Gamma_i^{\widehat{X}}\right)\subset\Gamma_i^X$.

%The homomorphism $a$ is trivial on the images of $\pi_1M$ and $\pi_1M^\tau$, thus $\widehat{X}$ contains $2m$ copies $M_1,\ldots,M_{2m}$
%of $M$ resp.\ $2m$ copies $M_1^\tau,\ldots,M_{2m}^\tau$
%of $M^\tau$. 
$\widehat{X}$ contains a copy of $\Sigma\times{\bf S}^1$ which is a $2m$-fold covering of $T^\tau$.
Let
$\widehat{M}=\pi^{-1}\left(M\right),
\widehat{M}^\tau=\pi^{-1}\left(M^\tau\right)$. Since $a\mid_{\pi_1M}$ is trivial, $\widehat{M}$ consists of $2m$ copies of $M$. If $\Sigma$ is separating, then $a\mid_{\pi_1M^\tau}$ is trivial, thus
also $\widehat{M}^\tau$
consists of $2m$ copies of $M^\tau$. If $\Sigma$ is non-separating,
then 
$a\mid_{\pi_1M^\tau}$ is surjective, thus
$\widehat{M}^\tau$ is a connected $2m$-fold covering of $M^\tau$.
%The Seifert-van Kampen theorem 
%implies that $\pi_*\left(\Gamma_i^{\widehat{X}}\right)$ is generated by $t^{2m}$ and elements of $\Gamma_i$.
%From \hyperref[elliptic]{Observation \ref*{elliptic}} 
%we obtain
%$\rho_X\left(t^{2m}\right)=A^{2m}={\bf 1}$, hence
%$\rho_{\widehat{X}}
%\left(\Gamma_i^{\widehat{X}}\right)=\rho_X\left(\Gamma_i^X\right)=
%\rho_X\left(\Gamma_i\right)\subset Fix\left(c_i\right)$. 

Consider the transfer map $tr:H_*\left(X,\partial X;{\bf Q}\right)
\rightarrow H_*\left(\widehat{X},\partial \widehat{X};{\bf Q}\right)$ of the finite covering $\widehat{X}\rightarrow X$. Application of the transfer map to
\hyperref[a]{Equation \ref*{a}} yields
\begin{equation}\label{tr} 
i_{\widehat{M}*}\left[\widehat{M},\partial \widehat{M}\right]-
i_{\widehat{M}^\tau *}\left[\widehat{M}^\tau,\partial \widehat{M}^\tau
\right]= i_{\Sigma\times{\bf S}^1 *}\left[\Sigma\times{\bf S}^1,\partial\Sigma\times{\bf S}^1\right]\end{equation}
in $H_3\left(\widehat{X},\partial \widehat{X};{\bf Q}\right)$, where $i_*$ denotes the respective inclusions
into $\widehat{X}$.

Composition of $\pi_*:\pi_1\widehat{X}\rightarrow \pi_1X$ with the representation $\rho_X:\pi_1X\rightarrow SL\left(2,{\bf C}\right)$
yields a representation $\rho_{\widehat{X}}:\Gamma^{\widehat{X}}\rightarrow 
SL\left(2,{\bf C}\right)$. 
By \hyperref[elliptic]{Observation \ref*{elliptic}} we have $A^m=\pm {\bf 1}$, thus $A^{2m}={\bf 1}$.
Recall that $\rho_X\left(t\right)=A$, thus $\rho_{\widehat{X}}\left(t^{2m}\right)=A^{2m}={\bf 1}$.
% sends the generator $t^{2m}$
% of $\pi_1{\bf S}^1\subset\pi_1\left(\Sigma\times{\bf S}^1\right)$ to $A^{2m}={\bf 1}$.

The Seifert-van Kampen Theorem
implies that $\pi_*\left(\Gamma_i^{\widehat{X}}\right)$ is generated by $t^{2m}$ and elements of $\Gamma_i$. Since
$\rho_X\left(t^{2m}\right)=A^{2m}={\bf 1}$ this implies $\rho_{\widehat{X}}
\left(\Gamma_i^{\widehat{X}}\right)=
%\rho_X\left(\Gamma_i^X\right)=
\rho_X\left(\Gamma_i\right)\subset Fix\left(c_i\right)$.
%% with 
%$\rho_{\widehat{X}}\left(\Gamma_i^{\widehat{X}}\right)\subset Fix\left(c_i\right)$. 
Thus we can 
extend $$B\rho_{\widehat{X}}:\left(B\Gamma^{\widehat{X}}\right)^{comp}\rightarrow BSL\left(2,{\bf C}\right)^{comp},$$ 
by mapping the cone point over $B\Gamma_i^{\widehat{X}}$ to $c_i$, such that $B\rho_{\widehat{X}}Bj=B\rho$ for the homomorphism of pairs $j:\left(\Gamma,\Gamma_i\right)\rightarrow \left(\Gamma^{\widehat{X}}, \Gamma^{\widehat{X}}_i\right)$ induced by the inclusion $i:\widehat{M}\rightarrow \widehat{X}$. 
%Let $\rho_1:=\rho_{\widehat{X}} \circ \hat{j}_{\wide}, \rho_l^\tau:=\rho_{\widehat{X}}\circ \hat{j}_{M_l^\tau}$ for the homomorphisms $\hat{j}_*$ induced on $\pi_1$ by the composition of $\hat{i}_*$ with the lift $M\rightarrow M_l$. Then we have that $B\rho_{\widehat{X}}B\hat{j}_{M_l}$ agrees with $B\rho_l$

We do not know whether $X$ and $\widehat{X}$ are aspherical or not, but we do have the classifying map $\Psi_{\widehat{X}}:\widehat{X}\rightarrow 
\mid B\Gamma^{\widehat{X}}\mid$, whose restriction to $\partial_i\widehat{X}$ is (upon a homotopy) the classifying map 
$\partial_i
\widehat{X}\rightarrow 
\mid B\Gamma^{\widehat{X}}_i\mid$ and which therefore extends to $$\Psi_{\widehat{X}}:DCone\left(\cup_{i=1}^s\partial_i\widehat{X}\rightarrow \widehat{X}\right)
\rightarrow \mid\left(B\Gamma^{\widehat{X}}\right)^{comp}\mid.$$ 

Let $P_1:\Sigma\times {\bf S}^1\rightarrow \Sigma$ be the projection to the first factor, then $$\left(\mid B\rho_{\widehat{X}}\mid\Psi_{\widehat{X}}i_\Sigma P_1\right)\sim\left(\mid B\rho_{\widehat{X}}\mid\Psi_{\widehat{X}} i_{\Sigma
\times{\bf S}^1}\right):DCone\left(\cup_j\partial_j\Sigma\times{\bf S}^1
\rightarrow\Sigma\times{\bf S}^1\right)\rightarrow \mid BSL\left(2,{\bf C}\right)^{comp}\mid$$ are homotopic
% the
%diagram
%$$\begin{xy}
%\xymatrix{ \Sigma\times {\bf S}^1\ar[r]^{P_1} \ar[d]^{p_m} & \Sigma\ar[d]^{i_{\Sigma}}\\
%T^\tau\ar[r]^{i_{T^\tau}}&X
%}
%\end{xy}$$ commutes 
since $\Sigma\times{\bf S}^1$ and $\partial_j \Sigma\times{\bf S}^1$ are aspherical and since the induced homomorphisms between fundamental groups agree because $\rho_{\widehat{X}}$ sends the generator $t^{2m}$
of $\pi_*\left(\pi_1{\bf S}^1\right)\subset\pi_*\left(\pi_1\left(\Sigma\times{\bf S}^1\right)\right)$ to $A^{2m}={\bf 1}$.

Since $$\left(\mid B\rho_{\widehat{X}}\mid\Psi_{\widehat{X}}i_{\Sigma\times{\bf S}^1}\right)_*=\left(\mid B\rho_{\widehat{X}}\mid\Psi_{\widehat{X}} i_{\Sigma} P_1\right)_*:H_3\left(\Sigma\times {\bf S}^1,\partial\Sigma\times{\bf S}^1\right)\rightarrow H_3\left(\mid BSL\left(2,{\bf C}\right)^{comp}\mid\right)$$
factors over $H_3\left(\Sigma,\partial\Sigma\right)=0$ we obtain $$\left(\mid B\rho_{\widehat{X}}\mid
\Psi_{\widehat{X}}i_{\Sigma\times{\bf S}^1}\right)_*\left[\Sigma\times {\bf S}^1,\partial\Sigma\times{\bf S}^1\right]=0.$$ 
%Application of $EM^{-1}$ yields then
%$$\left(\mid B\rho_{\widehat{X}}\mid\Psi_{\widehat{X}}\hat{i}_{\Sigma\times{\bf S}^1}\right)_*
%\left[\Sigma\times {\bf S}^1,\partial\Sigma\times{\bf S}^1\right]=0,$$
%EM^{-1}\left(i_{T^\tau}p_m\right)_*\left[\Sigma\times {\bf S}^1,\partial\Sigma\times{\bf S}^1\right]=0,$$
Thus
\hyperref[tr]{Equation \ref*{tr}} implies
\begin{equation}\label{d}
\left(\mid B\rho_{\widehat{X}}\mid\Psi_{\widehat{X}}i_{\widehat{M}}\right)_*
\left[\widehat{M},\partial\widehat{M}\right]=\left(
(\mid B\rho_{\widehat{X}}\mid\Psi_{\widehat{X}}
i_{\widehat{M}^\tau}\right)_*\left[\widehat{M}^\tau,\partial \widehat{M}^\tau\right]\in H_3\left(\mid B
%\Gamma^{\widehat{X}}
SL\left(2,{\bf C}\right)^{comp};{\bf Q}\mid\right).\end{equation}

%We have that $$\mid B\rho_{\widehat{X}}\mid\Psi_{\widehat{X}}i_{\widehat{M}}\sim\mid B\rho_{X}\mid\Psi_{X} i_{M}\pi:
%DCone\left(\cup_i\partial_i\widehat{M}
%\rightarrow\widehat{M}\right)\rightarrow \mid BSL\left(2,{\bf C}\right)^{comp}\mid$$ are homotopic because $\widehat{M}$ and $\partial_i\widehat{M}$ are aspherical and the induced homomorphims between fundamental groups agree.
%Therefore, \hyperref[d]{Equation \ref*{d}} together with $\pi_*\left[\widehat{M},\partial\widehat{M}\right]=m\left[M,\partial M\right]$ and 
%$\pi_*\left[\widehat{M}^\tau,\partial\widehat{M}^\tau\right]=m\left[M^\tau,\partial M^\tau\right]$ implies 

%Application of $\left(B\rho_{\widehat{X}}\right)_*$ yields, 
Let $\widehat{M}_1$ be the path-component of $\widehat{M}$ with $\hat{x}\in \widehat{M}_1$, let
$i_{\widehat{M}_1}:\widehat{M}_1\rightarrow \widehat{X}$ be the inclusion and $I_1:M\rightarrow\widehat{M}_1$ 
the homeomorphism inverse to $\pi\mid_{\widehat{M}_1}$. Let $j_1=\left(i_{\widehat{M}_1}I_1\right)_*:\Gamma\rightarrow \Gamma^{\widehat{X}}$ be the induced homomorphism, then we have $\rho=\rho_{\widehat{X}}j_1$ and therefore
%With $\hat{\rho}:=\rho_{\widehat{X}} \circ j_{\widehat{M}}, \hat{\rho}^\tau:=\rho_{\widehat{X}}\circ j_{\widehat{M}^\tau}$ for the homomorphisms $j_*$ induced on $\pi_1$ by $i_*$ we have 
$$\left(\mid B\rho_{\widehat{X}}\mid\Psi_{\widehat{X}}i_{\widehat{M}_1}I_1\right)_*=
\left(\mid B\rho_{\widehat{X}}\mid \mid Bj_1\mid \Psi_{M}\right)_*=
\left(\mid B\rho\mid\Psi_{M}\right)_*=\left(\mid B\rho\mid\right)_*\iota^{-1} EM^{-1}=\iota^{-1}\left(B\rho\right)_*EM^{-1}$$
(where $\iota:H_*\left(\mid .\mid\right)\rightarrow H_*^{simp}\left(.\right)$ denotes the isomorphism between singular and simplicial homology).

If $\widehat{M}_l$ is another path-component of $\widehat{M}$, 
then we have a deck transformation $\sigma_l:\widehat{X}\rightarrow\widehat{X}$ 
with $\sigma_l\left(\widehat{M}_1\right)=\widehat{M}_l$. We denote $I_l=\sigma_lI_1$ and obtain by the same computation 
$\left(\mid B\rho_{\widehat{X}}\mid\Psi_{\widehat{X}}i_{\widehat{M}_l}I_l\right)_*=\iota^{-1}\left(B\rho\right)_*EM^{-1}$.
%The composition 
%$\Psi_{\widehat{M}_l}\sigma_l$ is a classifying map for $\pi_1\widehat{M}_1$, 
%therefore $\Psi_{\widehat{M}_l}\sigma_l\sim 
%\Psi_{\widehat{M}_1}$ and thus $$\left(\mid B\rho_{\widehat{X}}
%\mid\Psi_{\widehat{X}}i_{\widehat{M}_l}\right)_*\left[\widehat{M}_l,\partial \widehat{M}_l\right]=\left(\mid B\rho_{\widehat{X}}\mid\Psi_{\widehat{X}}i_{\widehat{M}_1}\right)_*\left[\widehat{M}_1,\partial\widehat{M}_1\right]=
%\iota^{-1}\left(B\rho\right)_*EM^{-1}\left[\widehat{M},\partial\widehat{M}_1\right].$$
If  
{\em $\Sigma$ is separating}, then the same argument applies to the path-components of $\widehat{M}^\tau$: let $\widehat{M}_1^\tau$ be the path-component of $\widehat{M}^\tau$ with $\hat{x}\in\widehat{M}_1^\tau$, let
$I_1^\tau:M^\tau
\rightarrow\widehat{M}_1^\tau$ be 
the homeomorphism inverse to $\pi\mid_{\widehat{M}_1^\tau}$ and $I_l^\tau=\sigma_lI_1^\tau$, then we obtain $\left(\mid B\rho_{\widehat{X}}\mid\Psi_{\widehat{X}}i_{\widehat{M}_l^\tau}I_l^\tau\right)_*=\iota^{-1}\left(B\rho^\tau\right)_*EM^{-1}$.

%, and the analogous equality for $\widehat{M}^\tau$.
We have $\left[\widehat{M},\partial \widehat{M}\right]=\sum_{l=1}^{2m} I_{l*}\left[M,\partial M\right]$ and, if $\Sigma$ is separating, also
$\left[\widehat{M}^\tau,\partial \widehat{M}^\tau\right]=\sum_{l=1}^{2m} I_{l*}^\tau\left[M^\tau,\partial M^\tau\right]$.
% for the analogously defined $I_l^\tau$. 

Thus \hyperref[d]{Equation \ref*{d}} implies
$$2m \left(B\rho\right)_* EM^{-1}\left[M,\partial M\right]=
2m \left(B\rho^\tau\right)_* EM^{-1}\left[M^\tau,\partial M^\tau\right]$$
if $\Sigma$ is separating, 
respectively 
$$2m
\left(B\rho\right)_*
EM^{-1}\left[M,\partial M\right]=
\left(B\hat{\rho}^\tau\right)_* EM^{-1}\left[\widehat{M}^\tau,\partial\widehat{M}^\tau\right]$$
with $\hat{\rho}^\tau:=\rho_{\widehat{X}}\left(i_{\widehat{M}^\tau}\right)_*$ 
if $\Sigma$ is non-separating. (In the latter case $\widehat{M}^\tau$ was connected.)

%\begin{equation}\label{e}B\rho_*EM^{-1}\left[M,\partial M\right]=B\rho^\tau_*EM^{-1}\left[M^\tau,\partial M^\tau\right]\in H_3\left(BSL\left(2,{\bf C}\right)^{comp};{\bf Q}\right).\end{equation}

By Section 3.1 we have $$<\left[\overline{c\nu_3}\right], \left(B\rho\right)_*EM^{-1}\left[M,\partial M\right]>=vol\left(M\right),$$ 
 $$<\left[\overline{c\nu_3}\right], \left(B\rho^\tau\right)_*EM^{-1}\left[M^\tau,\partial M^\tau\right]>=
 vol\left(M^\tau\right),$$
and, if $\Sigma$ is non-separating,
 $$<\left[\overline{c\nu_3}\right], \left(B\hat{\rho}^\tau\right)_*EM^{-1}\left[\widehat{M}^\tau,\partial \widehat{M}^\tau\right]>=vol\left
(\widehat{M}^\tau\right)=2m vol\left(M^\tau\right)$$. 

%Thus $\left(B\rho\right)_*EM^{-1}\left[\widehat{M},\partial \widehat{M}\right]=\left(B\rho\right)_*EM^{-1}\left[\widehat{M}^\tau,\partial \widehat{M}^\tau\right]$ implies 
Hence $vol\left(M\right)=vol\left(M^\tau\right)$.

\end{pf}\\
\\
Remark: The assumption that $\Sigma$ is not a virtual fiber is needed "only" for the proof that $M^\tau$ is hyperbolic, that is for the application of \hyperref[discrete]{Proposition \ref*{discrete}}.

If $\Sigma$ is a virtual fiber then $M^\tau$ may or may not be hyperbolic. For example, if $M=T^\alpha$ is a mapping torus of $\alpha:\Sigma\rightarrow\Sigma$, then $M^\tau$ is the mapping torus of $\alpha\circ\tau:\Sigma
\rightarrow\Sigma$. By Thurstons hyperbolization theorem, the mapping torus $T^\alpha$ is hyperbolic if and only if $\alpha$ is pseudo-Anosov. However, if $\alpha$ is pseudo-Anosov and $\tau$ is of finite order, then $\alpha
\circ\tau$ may or may not be pseudo-Anosov, so there is no general statement whether $M^\tau$ is hyperbolic or not.

Examples for both phenomena can be easily found in $SL\left(2,{\bf Z}\right)$, the mapping class group of the once-punctured
torus. Here Anosov diffeomorphisms correspond to hyperbolic elements in
$SL\left(2,{\bf Z}\right)$ and finite order diffeomorphisms correspond to elliptic elements in $SL\left(2,{\bf Z}\right)$.  One easily finds hyperbolic elements $A_1,A_2$ and elliptic elements $B_1,B_2$ such that $A_1B_1$ is hyperbolic while $A_2B_2$ is not.

However, if $M^\tau$ happens to be hyperbolic, then the proof of \hyperref[cusped]{Theorem \ref*{cusped}} shows that $vol\left(M^\tau\right)=vol\left(M\right)$ even if $\Sigma$ was a virtual fiber.

%\noindent
%Thilo Kuessner\\
%Mathematisches Institut, Universit\"at M\"unster\\
%Einsteinstra\ss e 62\\
%D-48149 M\"unster\\
%Germany\\
%e-mail: kuessner@math.uni-muenster.de\\

\end{document}